\documentclass{article}
\usepackage{amssymb}
\usepackage{latexsym}
\newcommand{\N}{{\mathbb N}}

\newcommand{\Z}{{\mathbb Z}}

\newcommand{\bx}{\hfill{$\Box $ }}

\newtheorem{corollary}{Corollary}[section]
\newtheorem{proposition}{Proposition}[section]
\newtheorem{theorem}{Theorem}[section]
\title{Some Properties of Finite Rings}
\author{R.Coleman\\Laboratoire Jean Kuntzmann\\
Domaine Universitaire de Saint-Martin-d'H\`eres, France.}
\date{}
\begin{document}
\maketitle

\begin{abstract} A well-known theorem of Wedderburn asserts that a finite division ring is commutative. In a division ring the group of invertible elements is as large as possible. Here we will be particularly interested in the case where this group is as small as possible, namely reduced to $1$. We will show that, if this is the case, then the ring is boolean. Thus, here too, the ring is commutative.\\ 

\noindent Classification: 16K99, 16N99

\end{abstract}

\noindent {\bf Notation} We will write $|S|$ for the cardinal of a set $S$. If
$R$ is a ring, then we will denote the subset of its nonzero elements $R^*$ and
the subset of its invertible elements $R^{\times}$. If the ring is a division ring, then $R^*=R^{\times}$. Also, we will write $char (R)$ for the characteristic of a ring $R$.
 
\section{Boolean rings} A ring $R$ is boolean if all its elements are 
idempotent, i.e., $x^2=x$ for all $x\in R$. A simple example of a boolean ring 
is $\Z _2$. Products of boolean rings are also boolean, so we may construct a
large class of such rings. 

\begin{proposition} If $R$ is a boolean ring, then $char (R)=2$, $R$ is 
commutative and $R^{\times}=\{1\}$.
\end{proposition}

\noindent \textsc{proof} We have
$$
x+y = (x+y)^2 = x^2+xy+yx +y^2 = x+xy+yx +y,
$$
which implies that $xy+yx=0$. If we set $x=y=1$, then we obtain $char (R)=2$.

Now $char (R)=2$ implies that $xy+xy=0$. This, with the fact that $xy+yx=0$,
implies that $xy=yx$, i.e., $R$ is commutative.

Suppose now that $x\in R^{\times}$. Then, multiplying the expression $x^2=x$ by
$x^{-1}$, we obtain $x=1$. Thus $R^{\times}$ contains the unique element 1. This
finishes the proof.\bx\\

We should point out here that not all rings of characteristic 2 are Boolean. For
example, the ring ${\cal M}_2(\Z _2)$ of square $2\times 2$ matrices, with
coefficients in $\Z_2$, is not Boolean. The polynomial ring $\Z_2[X]$ is another  example.

\section{The group of invertible elements}

The invertible elements of a ring form a group with the multiplication of the 
ring. In this section we will consider some elementary properties of this group.

\begin{proposition}\label{prop.inv1} Let $R$ be a ring whose characteristic is not 2. If $x$ is invertible, then $-x\neq x$. It follows that, if $R$ is finite, then the sum of the elements of $R^{\times}$ is 0.
\end{proposition} 

\noindent \textsc{proof} Let $x\in R^{\times}$. Then
$$
x=-x \Longrightarrow  x+x=0 \Longrightarrow 1+1=0,
$$
a contradiction, because $char (R)\neq 2$. It follows that $x\neq -x$. 

Suppose now that $R$ is finite. If $x$ is invertible, then so is $-x$. Therefore
$R^{\times}$ is composed of pairs whose sum is 0. Thus the sum of the elements
of $R^{\times}$ is 0.\bx

\begin{corollary} If $char(R)\neq 2$ and $R$ is finite, then $|R^{\times}|$ is 
an even number.
\end{corollary}

\noindent {\bf Remark.} If $x\in R$ is not invertible, then we may have $x=-x$, 
even if the characteristic of the ring is not 2. For example, in $\Z _4$, which 
is of chacteristic 4, $2=-2$. In fact, more generally in $\Z _{2n}$, $n=-n$.\\

We may extend Proposition \ref{prop.inv1} to finite fields of more than two elements, even if the chacteristic is 2.
  
\begin{proposition} Let $R$ be a finite field. If $|R|>2$, then the sum of
the elements of $R^{\times}$ is 0. If $|R|=2$, then the sum is 1.
\end{proposition}

\noindent \textsc{proof} If $|R|=2$, then $R^{\times}$ contains the unique
element 1, hence the result.

Suppose now that $|R|=n>2$ and that $\alpha $ is a generator of the group
$R^{\times}$. Then 
$$
0 = 1-\alpha ^n = (1-\alpha )(1 + \alpha + \cdots + \alpha ^{n-1}).
$$
As $\alpha \neq 1$, we have 
$$
1 + \alpha + \cdots + \alpha ^{n-1} =0.
$$
However, this is the sum of the elements of $R^{\times}$. Hence the result.
\bx\\

\noindent {\bf Remark.} If a finite field $R$ is of characteristic $2$, then $R$ has $2^s$ elements for some $s\in \N^*$. Hence $|R^{\times}|$ is an odd number. This is not the case if the characteristic is an odd prime number.

\section{Matrix rings over finite fields}

In this section we consider the particular case of matrix rings over finite fields of characteristic 2. We will write ${\cal M}_n(R)$ for the set of $n\times n$ matrices with coordinates in $R$. With the usual operations of addition and multiplication of matrices, ${\cal M}_n(R)$ is a ring. If $n=1$, 
then ${\cal M}_n(R)$ is isomorphic to $R$. If $|R|=2$, then $R^{\times}$ contains the unique element 1, otherwise the sum of its elements is
0. Now let us consider the case where $n>1$. In this case ${\cal M}_n(R)$ is
noncommutative.

\begin{proposition} If $R$ is a finite field of characteristic 2 and $n\geq 2$, then the sum of the elements of ${\cal M}_n(R)^{\times}$ has an even number of
elements, whose sum is 0.
\end{proposition}

\noindent \textsc{proof} If $R$ is a finite field of characteristic $p$, then 
$$
|{\cal M}_n(R)^{\times}|= (p^n-p^{n-1})(p^n-p^{n-2})\cdots (p^n-1).
$$
A proof may be found in \cite{rotman}. It follows that in the case where $p=2$ and $n>1$, $|{\cal M}_n(R)^{\times}|$ is an even number.

Let $c$ be a nonzero vector in $R^2$. As $|R|=2^s$ for a certain $s\in \N^*$, there are $2^s$ multiples of $c$. Therefore there are $2^{2s}-2^s$ vectors which are not multiples of $c$. Thus there is an even number of matrices in 
${\cal M}_2(R)^{\times}$ having the first column $c^T$. A similar argument to that we have used for the first column shows that there is an even number of matrices in ${\cal M}_2(R)^{\times}$ having the same second
column. It follows that the sum of the elements in ${\cal M}_2(R)^{\times}$ is 0.

Now let us consider the case $n>2$. Let $c_1$ be a nonzero 
vector in $R^n$ and $c_2,\ldots ,c_n\in R^n$ be such that $c_1,\ldots ,c_n$ 
form an independant set. If we fix $c_1$ and permute the other elements, then 
we obtain another distinct ordered set. There are $(n-1)!$ such permutations. 
Thus there are $(n-1)!$ matrices in ${\cal M}_n(R)^{\times}$ having the same 
columns with the first column fixed. As 2 divides $(n-1)!$, there is an even 
number of matrices with the same first column.  The preceding argument applies to any column and so the sum of the matrices in ${\cal M}_n(R)^{\times}$ is 0.\bx

\section{The main theorem} 

Our aim in this section is to show that a finite ring $R$ in which the multiplicative group is as small as possible, i.e., $R^{\times}=\{1\}$, is a Boolean ring.\\

We will first give a brief review of Artin-Wedderburn theory. The subject is well-handled in various places. A good reference is \cite{ash}.\\ 

A division ring is a ring $R$ such that $R^{\times}=R^*$. A theorem of Wedderburn states that, if such a ring is finite, then it is commutative, i.e., a field. Proofs of this result may be found, for example,  in \cite{herstein} or \cite{lam}. It is natural to consider another 'extreme' case, namely where $R^{\times}=\{1\}$. Our aim in this section is to show that in this case $R$ is a boolean ring.\\

We will first recall some definitions and results from elementary ring theory. In a ring $R$ the intersection of the maximal left ideals is called the 
Jacobson radical of $R$ and usually noted $J(R)$. It turns out that $J(R)$ is also the intersection of all the maximal right ideals and so is an ideal. We may characterize elements of $J(R)$ in the following way: $a\in J(R)$ if and only if
$1-xay\in R^{\times}$ for all $x,y\in R$. By setting $x=y=1$, we see that if $R^{\times}=\{1\}$, then $J(R)=\{0\}$.\\

We now recall that a ring is artinian if any descending chain of ideals is stationary after a finite number of ideals. If we replace ideals by left (resp. right) ideals in the definition, then we obtain the definition of a left (resp. right) artinian ring. Clearly a finite ring is artinian, as well as being left and right artinian.\\

We now come to the the notion of semi-simplicity and to Wedderburn's structure theorem. In the definitions we will use left ideals; however, we could replace these by right, or two-sided, ideals. To be brief, we will use the term ideal for left ideal. We say that an ideal $I$ is simple if $I\neq \{0\}$ and the only sub-ideals included in $I$ are $I$ itself and $\{0\}$. A ring is semi-simple if it is a product of simple ideals. It should be noticed that a ring is semi-simple if and only if it is a direct sum of simple ideals. A fundamental result in the theory of semi-simple rings is the Wedderburn structure theorem, namely

\begin{theorem} If $R$ is a semi-simple ring, then there are division rings $D_1,\ldots ,D_t$ and positive integers $n_1, \ldots ,n_t$ such that 
$$
R\simeq {\cal M}_{n_1}(D_1)\oplus \cdots \oplus {\cal M}_{n_t}(D_t).
$$
\end{theorem}

If $R$ is an artinian ring and $I$ an ideal in $R$, then the quotient ring $R/I$ is also artinian. In the case where $I=J(R)$ the ring $R/I$ is semi-simple. Of course, if $J(R)=\{0\}$, then $R$ is semi-simple. We thus have the following result:

\begin{theorem} If $R$ is a finite ring and $R^{\times}=\{1\}$, then $R$ is semi-simple.
\end{theorem}

As a corollary we have the principle result of this paper, namely

\begin{theorem} If $R$ is a finite ring and $R^{\times}=\{1\}$, then $R$ is a boolean ring.
\end{theorem} 

\noindent \textsc{proof} As $R$ is semi-simple, we know that there are division rings $D_1,\ldots ,D_t$ and positive integers $n_1,\ldots ,n_t$ such that
$$
R\simeq {\cal M}_{n_1}(D_1)\oplus \cdots \oplus {\cal M}_{n_t} (D_t).
$$
As $R$ is finite so are the division rings. From Wedderburn's Little Theorem these rings are fields. Given that $|R^{\times}|=1$, it must be so that $D_i=\Z_2$ and $n_i=1$ for all $i$. However, ${\cal M}_1(\Z_2)\simeq \Z_2$ and so $R$ is isomophic to a direct sum of boolean rings and hence is a boolean ring.\bx\\

\noindent {\bf Remark.} As a boolean ring is commutative, if a finite ring $R$ is such that $R^{\times}=\{1\}$, then $R$ is commutative. Thus when $R^{\times}$ is as small as possible or as large as possible $R$ is commutative. 

\section{A final comment}

From what we have seen above we might be tempted to think that the sum of the invertible elements in a finite ring is either 0 or 1. However, this is not the case. We only need to consider the subring $A$ of upper triangular matrices in ${\cal M}_2(\Z _2)$. Here $A^{\times}$ is composed of the identity matrix $I_2$ and the matrix $M=(m_{ij})$, with $m_{21}=0$ and $m_{ij}=1$ otherwise. The sum of these matrices is the matrix $N=(n_{ij})$, with $n_{12}=1$ and $n_{ij}=0$ otherwise.\\ 

It should be noted that the sum of the invertible upper triangular matrices in ${\cal M}_n(\Z _2)$, with $n\geq 3$, is 0. There are $2^{\frac{(n-1)n}{2}}$ such matrices and, for $i<j$, in the position $ij$ half of them have a 0 and half of them a 1. As the number $2^{\frac{(n-1)n}{2}}$ is divisible by 4, there is an even number of matrices with 1 in any position $ij$, with $i<j$, and it follows that the sum of the matrices which interest us is $0$.


\begin{thebibliography}{10}
\bibitem{ash}
Ash R.,
{\it Basic Abstract Algebra},
Dover 2006
\bibitem{cohn}
Cohn P.M.,
{\it Introduction to the Theory of Rings},
Springer-Verlag 2000
\bibitem{herstein}
Herstein I.N.,
{\it Topics in Algebra},
John Wiley and Sons 1975
\bibitem{lam}
Lam T.Y.,
{\it A First Course in Noncommutative Rings},
Springer-Verlag 2001
\bibitem{rotman}
Rotman J.L.,
{\it An Introduction to the Theory of Groups},
Springer-Verlag 1999



\end{thebibliography}
\end{document}